\documentclass[twoside]{article}
\input amssym.def
\input amssym.tex
\usepackage{bbm}
\usepackage{graphicx}
\usepackage{color}
\usepackage{epsfig}
\usepackage{psfrag}

\skip\footins 22pt plus 2pt minus 2pt
\newdimen\jot
\def\mqth{\mathsurround=0pt }
\dimendef\dimenq=0
\def\openup{\afterassignment\qpenup\dimenq=}
\def\qpenup{\advance\lineskip\dimenq
  \advance\baselineskip\dimenq \advance\lineskiplimit\dimenq}
\def\eqalign#1{\,\vcenter{\openup1\jot \mqth
  \ialign{\strut\hfil$\displaystyle{##}$&$\displaystyle{{}##}$\hfil
  \crcr#1\crcr}}\,}
\newif\ifdtqp
\def\displqy{\global\dtqptrue \openup1\jot \mqth
  \everycr{\noalign{\ifdtqp \global\dtqpfalse
     \vskip-\lineskiplimit \vskip\normallineskiplimit
     \else \penalty\interdisplaylinepenalty \fi}}}
\def\displaylines#1{\displqy
  \halign{\hbox to\displaywidth{$\hfil\displaystyle##\hfil$}\crcr
  #1\crcr}}
\newskip\centerinq \centerinq=0pt plus 1000pt minus 1000pt

\newcommand{\ph}{\varphi}
\newcommand{\dast}{{\displaystyle{\ast}}}
\newcommand{\reals}{{\mathbbm{R}}}

\newcommand{\Lk}{\mathop{\rm Lk}\nolimits}
\newcommand{\area}{\mathop{\rm area}\nolimits}

\newcommand{\He}{\mathop{\raise.02ex\hbox{\rm H}}\nolimits}
\newcommand{\K}{\mathop{\raise.02ex\hbox{\rm K}}\nolimits}
\newcommand{\G}{\mathop{\raise.02ex\hbox{\rm G}}\nolimits}

\newcommand{\vol}{\mathop{\rm vol}\nolimits}

\newcommand{\Int}{\mathop{\rm Int}\nolimits}
\newcommand{\x}{{\bf x}}
\newcommand{\y}{{\bf y}}

\renewcommand{\t}{{\bf t}}
\newcommand{\s}{{\bf s}}

\newcommand{\A}{{\bf A}}
\newcommand{\p}{{\bf p}}

\begin{document}
\thispagestyle{empty}
\large

\centerline{\LARGE{\bf Linking integrals in the $n$-sphere}}

\medskip

\centerline{\sc Dennis DeTurck and Herman Gluck}

\bigskip

\bigskip

\small
\centerline{\bf ABSTRACT\rm}
\begin{quotation}
\noindent Let $K$ and $L$ be disjoint closed oriented submanifolds of the $n$-sphere $S^n$, with dimensions adding up to $n-1$. We define a map from their join $K*L$ to $S^n$ whose degree up to sign equals their linking number, and then use this
to find the desired linking integral.

\medskip

\end{quotation}

\large

\addtolength{\baselineskip}{2pt}

\bigskip

Here is our main result.

\noindent\bf Theorem\it.\ \ Let $K^k$ and $L^\ell$ be disjoint closed oriented smooth
submanifolds of $S^n$ with $k+\ell=n-1$. Then their linking number is given by the integral
$$\Lk(K^k,L^\ell)={1\over\vol S^n}\int_{K\times L} {\ph_{k,\ell}(\alpha)
\over\sin^n\alpha}\, [\x,d\x,\y,d\y]$$
where 
$$\ph_{k,\ell}(\alpha)=\int_{\beta=\alpha}^\pi \sin^k(\beta-\alpha)\sin^\ell\beta\,
d\beta,$$
and $\alpha=\alpha(\x,\y)$ is the distance in $S^n$ between $\x$ and $\y$.\rm

\medskip

In this formula, $\x\in K^k$ and $\y\in L^\ell$. We explain the notation 
$[\x,d\x,\y,d\y]$ in \S~4, but simply mention here that it is an $(n+1)\times
(n+1)$ determininant in which $\x$ occupies one column, $d\x$ occupies $k$ columns,
$\y$ occupies one column, and $d\y$ occupies $\ell$ columns.

The above integral is geometrically meaningful in the sense that its integrand is invariant under orientation-preserving isometries of $S^n$.

To prove this theorem, we consider the \it join \rm $K*L$ of the manifolds $K$ and $L$, which is obtained from the product $K\times L\times [0,1]$ by collapsing $K\times L\times 0$ to $K$ and $K\times L\times 1$ to $L$. Since each point $\x\in K$ is distinct from each point $\y\in L$, it follows that the points $\x$ and $-\y$ are not antipodal
in $S^n$, and hence can be connected there by a unique shortest geodesic arc. 

We define a map $f\colon K*L\to S^n$ by sending the line segment $\{(\x,\y,u)\,|\,0\le u\le 1\}$ connecting $\x$ and $\y$ in $K*L$ proportionally to the geodesic arc connecting $\x$ and $-\y$ in $S^n$. The degree of this map $f$ is, up to sign, the linking number $\Lk(K,L)$.

To evaluate this degree, and hence the linking number, we take the volume form on $S^n$, pull it back via $f$ to an $n$-form on $K*L$, partially integrate this $n$-form
along the line segments $\{(\x,\y,u)\,|\,0\le u\le 1\}$, and obtain the formula in the theorem above.

The anti-commutation rule 
$$\Lk(K^k,L^\ell)=(-1)^{(k+1)(\ell+1)}\Lk(L^\ell,K^k)$$
is reflected in the form of the linking integral, since
$$[\y,d\y,\x,d\x]=(-1)^{(k+1)(\ell+1)}[\x,d\x,\y,d\y],$$
as a result of interchanging $(k+1)$ columns with $(\ell+1)$ columns in these determinants, while
$$\ph_{k,\ell}(\alpha)=\ph_{\ell,k}(\alpha).$$

Suppose now that the submanifolds $K$ and $L$ of $S^n$ are disjoint, not only from one another, but also each from the antipodal image of the other. We comment in the next section on the origin of this hypothesis. In such a case we obtain the following

\noindent\bf Corollary\it.\ \ Let $K^k$ and $L^\ell$ be closed smooth oriented submanifolds of $S^n$ with \mbox{$k+\ell=n-1$}, and with $K$ disjoint from both $L$ and its antipodal image $-L$. Then
$$\Lk(K^k,L^\ell)+(-1)^n\Lk(K^k,-L^\ell)={(-1)^k\over \vol S^n}\int_{K\times L}
{\sin^k*\sin^\ell(\alpha)\over\sin^n\alpha}\,[\x,d\x,\y,d\y]$$ where the convolution $\sin^k*\sin^\ell$ is defined by 
$$\sin^k*\sin^\ell(\alpha) = \int_{\beta=0}^\pi\sin^k(\alpha-\beta)\sin^\ell\beta\ d\beta.$$\rm

In particular, if $K$ and $L$ are disjoint submanifolds lying in some open hemisphere of $S^n$, then $-L$ lies in the complementary hemisphere, and hence $\Lk(K,-L)=0$. In such a case we get
$$\Lk(K^k,L^\ell)={(-1)^k\over \vol S^n}\int_{K\times L}
{\sin^k*\sin^\ell(\alpha)\over\sin^n\alpha}\,[\x,d\x,\y,d\y].$$

All the linking integrals presented above were obtained at the same time by Clayton Shonkwiler and David Shea Vela-Vick (2008), using different methods, as a special case of their more general 
higher-dimensional linking integrals. We very much appreciate their help with the content and figures of the present paper.

\vfill
\eject

\noindent\bf \Large{1.\ Background}\rm
\large
\addtolength{\baselineskip}{2pt}

In a half-page paper dated January 22, 1833, Carl Friedrich Gauss gave without proof an integral formula for the linking number of two disjoint closed curves $K=\{\x(s)\}$ and $L=\{\y(t)\}$ in Euclidean 3-space,
$$\Lk(K,L)={1\over 4\pi}\int_{K\times L}{d\x\over ds}\times{d\y\over dt}\cdot {\x-\y\over|\x-\y|^3}\,ds\,dt.$$
The correspondence $(\x,\y)\to(\x-\y)/|\x-\y|$ defines a map from the torus $K\times L$ to the unit 2-sphere $S^2$ in $\reals^3$, whose degree up to sign is the linking number $\Lk(K,L)$. The area 2-form $\omega$ on $S^2$, when pulled back via this map to $K\times L$, integrated there and divided by $4\pi$, gives Gauss's formula. The mathematical historian Moritz Epple (1998) believes that this was the argument Gauss had in mind when he wrote the above formula.

This degree-of-map derivation of Gauss's linking integral works in Euclidean 3-space because the set of ordered pairs of distinct points in $\reals^3$ deformation retracts to a 
\mbox{2-sphere}. But it fails on the 3-sphere because the set of ordered pairs of distinct points in $S^3$ deformation retracts to a 3-sphere, and therefore all maps to it from the torus $K\times L$ are homotopically trivial.

Gauss undoubtedly knew another proof of his integral formula. Run a steady current through the first loop, and calculate the circulation of the resulting magnetic field around the second loop. By Amp\`ere's Law, this circulation is equal to the total current ``enclosed'' by the second loop, which means the current flowing along the first loop multiplied by the linking number of the two loops. Then the Biot-Savart forumla (1820) for the magnetic field leads directly to Gauss's linking integral.

To extend this line of thought, we developed in our (2005) paper a steady-state version of classical electrodynamics on the 3-sphere and in hyperbolic 3-space, and used this to obtain explicit integral formulas for the linking number of two disjoint curves in these spaces.

Greg Kuperberg, both in private correspondence and in his (2006) paper, independently used arguments completely different from ours to derive an equivalent linking integral on the 3-sphere. His arguments apply just as well on the $n$-sphere, and can be used to give another proof of our main theorem and its corollary.

Returning to the degree-of-map derivation of Gauss's linking integral, and aware that it does not extend to the 3-sphere, one can seek limited circumstances under which such an extension is possible, as follows.

Suppose that $K$ and $L$ are two smooth closed curves in $S^3$, disjoint from one another and also each from the antipodal image of the other.

Then, following Gauss, we have a natural map of the torus $K\times L$ into the set $P$ of pairs $(\x,\y)\in S^3\times S^3$ such that $\x\ne\y$ and $\x\ne-\y$. This set $P$ deformation retracts to the subset $P_0$ of orthogonal pairs $(\x,\y)\in S^3\times S^3$, which is just a copy of the unit tangent bundle $US^3$.

With that in mind, let $\omega$ be the $SO(4)$-invariant 2-form on $US^3$ which restricts to the area form on each fibre.

Let $f\colon K\times L\to US^3$ be the natural map of $K\times L$ into $P$, followed by its deformation retraction to $P_0=US^3$.

Then the pullback of $\omega$ via $f$ to $K\times L$ gives a 2-form which, when integrated over this torus and divided by $4\pi$, can be shown to be equal to the difference of the linking numbers $\Lk(K,L)-\Lk(K,-L)$. 

Shonkwiler and Vela-Vick 
(unpublished) have demonstrated that the same argument works on all odd-dimensional spheres (but not on the even-dimensional ones), and the result is the formula given in the corollary to our main theorem. Furthermore, they have shown how to deduce the formula for 
$S^n$ from that on $S^{n+1}$ to cover the even-dimensional case as well.

\bigskip

\noindent\bf \Large{2.\ Orientations and signs of linking numbers.}\rm
\large
\addtolength{\baselineskip}{2pt}

We discuss intersection and linking below in the combinatorial rather than the smooth setting, because null-homologous submanifolds bound chains, but not necessarily submanifolds, of one dimension higher. The definitions can be transplanted to the smooth category via smooth triangulations.

Let $M^n$ be an oriented $n$-dimensional combinatorial manifold, and let $P^p$ and $Q^q$ be subchains with $p+q=n$. If we assume that $P^p$ is a subcomplex of a given triangulation of $M^n$, and that $Q^q$ is a subcomplex of the Poincar\'e dual cell complex, then it follows that $P^p$ and $Q^q$ meet transversally in finitely many points, and we denote their intersection number by $\Int(P^p,Q^q)$.

Let $K^k$ and $L^\ell$ be two disjoint cycles in $M^n$ with $k+\ell=n-1$. Assume that $K$ and $L$ are null-homologous in $M$, with $K^k=\partial P^{k+1}$ and $L^\ell=\partial Q^{\ell+1}$, and that $K^k$ and $Q^{\ell+1}$ meet transversally as above, and likewise for $P^{k+1}$ and $L^\ell$. Then we define the \it linking number \rm of $K$ and $L$ to be 
$$\Lk(K^k,L^\ell)=\Int(K^k,Q^{\ell+1}).\eqno(2.1)$$
This definition can be seen to be independent of the choice of $Q^{\ell+1}$, thanks to the fact that $K$ is null-homologous in $M^n$. 

We note the inherent asymmetry in this definition, since one could equally well define the linking number of $K^k$ and $L^\ell$ to be the intersection number of $P^{k+1}$ and $L^\ell$. The two definitions can be shown to differ by the sign $(-1)^{k+1}$.

One can verify the anti-commutation rule,
$$\Lk(K^k,L^\ell)=(-1)^{(k+1)(\ell+1)}\Lk(L^\ell,K^k).\eqno(2.2)$$

We apply the definition (2.1) of linking number to the simplest possible example.

Orient Euclidean space $\reals^{n+1}$, and then orient the unit sphere $S^n\subset\reals^{n+1}$ so that a basis $\A_1,\ldots,\A_n$ for its tangent space at the point $\p$ is positive if and only if the basis $\p,\A_1,\ldots,\A_n$ for $\reals^{n+1}$ is positive. 

Assume that $\reals^{n+1}=\reals^{k+1}\oplus\reals^{\ell+1}$ is an orthogonal direct sum. 

Let the positive $x_1,\ldots,x_{k+1}$ axes in that order set the orientation for $\reals^{k+1}$, and orient its unit sphere $S^k$ as above. Likewise let the positive $y_1,\ldots,y_{\ell+1}$ axes in that order set the orientation for $\reals^{\ell+1}$, and orient its unit sphere $S^\ell$ as above. Let the positive $x_1\,\ldots,x_{k+1},y_1,\ldots,y_{\ell+1}$ axes in that order set the orientation for $\reals^{n+1}$, and then orient its unit sphere $S^n$ as above. 

Now $S^k$ and $S^\ell$ are oriented great subspheres of the oriented $S^n$, which link once geometrically. Using the definition (2.1), one checks easily that 
$$\Lk(S^k,S^\ell)=+1.\eqno(2.3)$$

\bigskip

\noindent\bf \Large{3.\ The degree of $f\colon K*L\to S^n$ and the linking number
$\Lk(K,L)$.}\rm
\large
\addtolength{\baselineskip}{2pt}

Let $K^k$ and $L^\ell$ be disjoint closed oriented smooth
submanifolds of $S^n$ with \mbox{$k+\ell=n-1$}. 
Let $f\colon K*L\to S^n$ be the map, defined in the introduction, which sends the line segment $\{(\x,\y,u)\,|\,0\le u\le 1\}$ connecting $\x$ and $\y$ in $K*L$ proportionally to the geodesic arc connecting $\x$ and $-\y$ in $S^n$.

A little elementary geometry yields the explicit formula
$$f(\x,\y,u)=\x\cos(u(\pi-\alpha))-{\y-\x\cos\alpha\over\sin\alpha}\sin(u(\pi-\alpha)),\eqno(3.1)$$
where $\alpha=\alpha(\x,\y)$ is the geodesic distance on $S^n$ between $\x$ and $\y$.

\begin{center}
	\begin{figure}
		\centering
	\psfrag{a}{\Large{$y$}}
	\psfrag{b}{\Large{$\alpha$}}
	\psfrag{c}{\Large{$x$}}
	\psfrag{e}{\Large{$\pi - \alpha$}}
	\psfrag{f}{\Large{$-y$}}
	\includegraphics[width=4in]{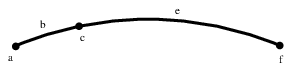}\\
	\large{Figure 1}
	\end{figure}
\end{center}

If $K^k$ and $L^\ell$ are closed oriented smooth manifolds, then the join $K^k*L^\ell$, although not in general itself a manifold, resembles an oriented $n$-manifold in the sense that it contains an $n$-cell as a dense open subset, and its top homology group is isomorphic to the integers. Orienting the join has the effect of choosing one or the other generator for its top homology, and we do this by orienting the product $K^k\times L^\ell\times[0,1]$ in the usual way, and then let this induce an orientation on $K^k* L^\ell$ under the collapsing map.

Now it makes sense to speak of the degree of the map $f\colon K^k*L^\ell\to S^n$ defined above, and we want to compare this degree to the linking number $\Lk(K^k,L^\ell)$. 

We try this out on the simplest possible example.

Let $K^k=S^k$ be the unit sphere in the $\reals^{k+1}$ with coordinates $x_1,\ldots,x_{k+1}$, let $L^\ell=S^\ell$ be the unit sphere in the $\reals^{\ell+1}$ with coordinates $y_1,\ldots,y_{\ell+1}$, and let $S^n$ be the unit sphere in the $\reals^{n+1}$ with coordinates $x_1,\ldots,x_{k+1},y_1,\ldots,y_{\ell+1}$. Assigning orientations as in the previous section, we noted there that $\Lk(S^k,S^\ell)=+1$.

The join of a $k$-sphere and an $\ell$-sphere is a sphere of dimension $k+\ell+1=n$, and it is easy to check that the natural map $S^k*S^\ell\to S^n$ has degree $+1$, \it provided \rm we give $S^k*S^\ell$ the orientation induced from that on the product $S^k\times[0,1]\times S^\ell$.

Our map $f\colon S^k*S^\ell\to S^n$ differs from this natural map in two ways. First, we orient $S^k*S^\ell$ via the product $S^k\times S^\ell\times[0,1]$ instead of $S^k\times[0,1]\times S^\ell$, which gives a sign change of $(-1)^\ell$ in the orientation. Second, we take the line segment 
$\{(\x,\y,u)\,|\,0\le u\le 1\}$ to the geodesic arc (quarter
circle in this case) connecting $\x$ to $-\y$, rather than connecting $\x$ to $\y$. 
The antipodal map on $S^\ell$ has degree $(-1)^{\ell+1}$, and so we get exactly this sign change in the degree.

The net effect of these two sign changes is simply multiplication by $-1$, and so we 
conclude that
$$\deg(f\colon S^k*S^\ell\to S^n)=-1=-\Lk(S^k,S^\ell).\eqno(3.2)$$
We will see in a moment the universal character of this example.

\eject

\noindent\bf Proposition 3.3\it.\ \ Let $K^k$ and $L^\ell$ be disjoint closed oriented smooth submanifolds of $S^n$ with $k+\ell=n-1$, and let $f\colon K^k*L^\ell\to\ S^n$ be the map described above. Then
$$\deg f = -\Lk(K^k,L^\ell).$$
\rm

\noindent\it Proof\rm.\ \ The simple idea of the proof is indicated in Figure 2, which shows 
a 2-component link in the 3-sphere. One component is a trefoil knot $K$ and the other is a circle $L$, with linking number $\Lk(K,L)=2$. Also shown in the figure is a Seifert surface $S$ (a punctured torus) bounded by $K$, and pierced twice transversally by $L$. 

\begin{center}
	\begin{figure}[h]
		\centering
	\psfrag{K}{\Large{\textcolor{red}{$K$}}}
	\psfrag{K1}{\Large{\textcolor{red}{$K'$}}}
	\psfrag{S}{\huge{\textcolor{red}{$S$}}}
	\psfrag{L}{\Large{\textcolor{blue}{$L$}}}
	\includegraphics[width=4in]{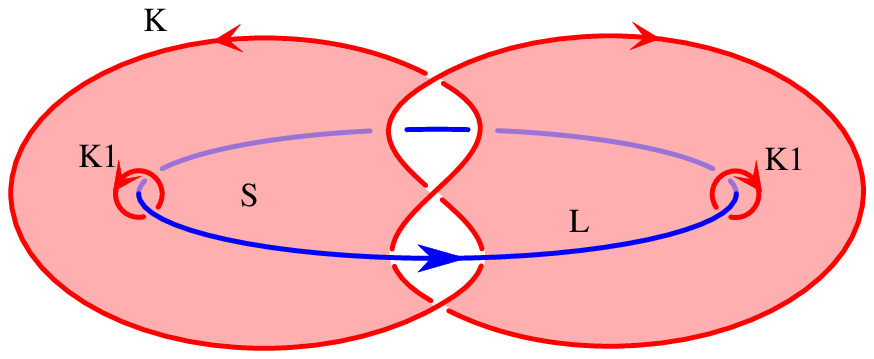}\\
	\large{Figure 2}
	\end{figure}
\end{center}

\vspace{-.3in}
Let $K'$ denote the union of the two small circles on $S$ which surround the punctures, let $D'$ denote the union of the two small disks on $S$ which they bound, and let $S'$ denote the remaining portion of $S$ bounded by $K$ and $K'$.

We have $\Lk(K,L)=\Lk(K',L)$ since $S\cap L=D'\cap L$.

Also, the maps $f\colon K*L\to S^3$ and $f'\colon K'*L\to S^3$ have the same degree, since they have a common extension to a map $F\colon S'*L\to S^3$.

Keeping $K'$ as is, if $L$ were not already a circle, we could repeat this construction by using a Seifert surface bounded by $L$ to replace it with a union $L'$ of small circles such that $\Lk(K',L)=\Lk(K',L')$ and such that the maps $f'\colon K'*L\to S^3$ and $f''\colon K'*L'\to S^3$ have the same degree.

The higher-dimensional version of this argument replaces $K$ by a union $K'$ of small round $k$-spheres, and $L$ by a union $L'$ of small round $\ell$-spheres, so that
$$\Lk(K,L)=\Lk(K',L')$$
and 
$$\deg(f\colon K*L\to S^n)=\deg(f''\colon K'*L'\to S^n).$$
Now our proposition follows from the special case given in (3.2), which is clearly just as valid for small round subspheres as for great ones. 

\vfill
\eject

\noindent\bf \Large{4.\ Proof of the main theorem.}\rm
\large
\addtolength{\baselineskip}{2pt}

Let $K^k$ and $L^\ell$ be disjoint closed oriented smooth
submanifolds of $S^n$ with \mbox{$k+\ell=n-1$}. 
Let $f\colon K*L\to S^n$ be the map defined in the introduction and given explicitly in (3.1). We saw in Proposition 3.3 that $\Lk(K^k,L^\ell)=-\deg f$, so our task is to find a good integral formula for the degree of $f$.

To do this, we start with the volume form $\omega$ on $S^n$, pull it back via $f$ to an $n$-form on $K*L$, partially integrate this $n$-form along the line segments
$\{(\x,\y,u)\,|\,0\le u\le 1\}$ so as to leave us with an integral over $K\times L$, and finally divide this integral by the volume of $S^n$ to get the degree of $f$.

Let $\s=(s_1,\ldots,s_k)$ be local coordinates on $K^k$, whose order gives the orientation there, and likewise for $\t=(t_1,\ldots,t_\ell)$ on $L^\ell$. 

The volume form $\omega$ on $S^n$ is given by 
$$\omega_{\p}(\A_1,\ldots,\A_n)=\det(\p,\A_1,\ldots,\A_n),\eqno(4.1)$$
according to our conventions about orientation in section 2.

Then the degree of $f$ is given by
$$\eqalign{
\deg f&= {1\over \vol S^n}\int_{K*L} f^\dast\omega \cr
&={1\over \vol S^n}\int_{K*L} f^\dast\omega\left({\partial\over\partial s_1},
\ldots,{\partial\over\partial s_k},{\partial\over\partial t_1},\ldots, {\partial\over \partial t_\ell},{\partial\over\partial u}\right)\,ds_1\cdots ds_k\,dt_1\cdots dt_\ell\,du.}\eqno(4.2)$$
We postpone integrating, and pay attention to the integrand,
$$\eqalign{
f^\dast\omega
&\left( 
{\partial\over\partial s_1},
\ldots,
{\partial\over\partial s_k},
{\partial\over\partial t_1},\ldots, 
{\partial\over \partial t_\ell},
{\partial\over\partial u}
\right)\cr
&=\omega_f\left({\partial f\over\partial s_1},
\ldots,{\partial f\over\partial s_k},{\partial f\over\partial t_1},\ldots, 
{\partial f\over \partial t_\ell},{\partial f\over \partial u}\right) \cr
&= \det\left(f,{\partial f\over\partial s_1},
\ldots,{\partial f\over\partial s_k},{\partial f\over\partial t_1},\ldots, 
{\partial f\over \partial t_\ell},{\partial f\over \partial u}\right) \cr
&= (-1)^{n-1}\det\left(f,{\partial f\over \partial u},{\partial f\over\partial s_1},
\ldots,{\partial f\over\partial s_k},{\partial f\over\partial t_1},\ldots 
{\partial f\over \partial t_\ell}\right).}\eqno(4.3)$$
We pay the price $(-1)^{n-1}$ to move the $\partial f/\partial u$ column adjacent to the $f$ column because it will be computationally convenient when we assess their joint contribution to the value of this determinant.

To that end, we take formula (3.1) for $f(\x,\y,u)$, differentiate it with respect to $u$, and get
$${\partial f\over\partial u}=-(\pi-\alpha)\left(\sin(u(\pi-\alpha))\x(\s)+
\cos(u(\pi-\alpha)){\y(\t)-\cos\alpha\,\x(\s)\over\sin\alpha}\right).\eqno(4.4)$$
We then compute that
$$f\wedge {\partial f\over \partial u}=-{\pi-\alpha\over\sin\alpha}\,\x(\s)\wedge \y(\t).\eqno(4.5)$$
The advantage of this preliminary computation is that when we compute $\partial f/\partial s_i$ and $\partial f/\partial t_j$, we can discard terms containing either $\x(\s)$ or $\y(\t)$, since they will disappear in the calculation of the determinant in (4.3) because of its alternating character. We signal this discard below with the symbol $\sim$. 

With this in mind, we compute that
$$
{\partial f\over \partial s_i}\sim \left(\cos(u(\pi-\alpha))+{\cos\alpha\over
\sin\alpha}\sin(u(\pi-\alpha))\right)\,{\partial\x\over \partial s_i}
={A\over \sin\alpha}\,{\partial \x\over\partial s_i}\eqno(4.6)$$
and
$${\partial f\over \partial t_j}\sim -\left({\sin(u(\pi-\alpha))\over
\sin\alpha}\right)\,{\partial\y\over \partial t_j}
=-{B\over \sin\alpha}\,{\partial \y\over\partial t_j},\eqno(4.7)$$
where we introduce for convenience the abbreviations
$$A=\sin\alpha\cos(u(\pi-\alpha))+\cos\alpha\sin(u(\pi-\alpha))\qquad\mbox{\rm and}
\qquad B=\sin(u(\pi-\alpha)).$$
Then, referring back to the last line of (4.3), we use (4.5), (4.6) and (4.7) to write
$$\eqalign{
\det &\left( f,{\partial f\over \partial u},
{\partial f\over\partial s_1},
\ldots,{\partial f\over\partial s_k},
{\partial f\over\partial t_1},\ldots, 
{\partial f\over \partial t_\ell}
\right)\cr
&=-{\pi-\alpha\over\sin\alpha}\left({A\over\sin\alpha}\right)^k
\left({-B\over\sin\alpha}\right)^\ell\det\left(\x,\y,{\partial\x\over \partial s_1},\ldots,
{\partial\x\over\partial s_k},{\partial\y\over\partial t_1},\ldots,
{\partial\y\over\partial t_\ell}\right)\cr
&=(-1)^n\,{\pi-\alpha\over\sin^n\alpha}\,A^kB^\ell\,\det\left(\x,
{\partial\x\over \partial s_1},\ldots,
{\partial\x\over\partial s_k},\y,{\partial\y\over\partial t_1},\ldots,
{\partial\y\over\partial t_\ell}\right),}\eqno(4.8)$$
since the cost of moving the $\y$ column to its new location in the determinant is $(-1)^k$.

Putting this all together, we have by (4.2), (4.3) and (4.8) that
$$\eqalign{
\deg &f\cr 
&={1\over\vol S^n}\int_{K*L}(-1)^{n-1}\det
\left( f,{\partial f\over \partial u},
{\partial f\over \partial s_1},\ldots,
{\partial f\over\partial s_k},
{\partial f\over\partial t_1},\ldots,
{\partial f\over\partial t_\ell}
\right)ds_1\cdots ds_k\, dt_1
\cdots dt_\ell\,du\cr
&={-1\over \vol S^n}\int_{K*L}{\pi-\alpha\over\sin^n\alpha}A^kB^\ell
\det\left(\x,{\partial\x\over\partial s_1},\ldots,{\partial\x\over\partial s_k},\y,
{\partial\y\over\partial t_1},\ldots,{\partial\y\over\partial t_\ell}\right)
ds_1\cdots ds_k\,dt_1\cdots dt_\ell\,du.}$$

As promised, we now partially integrate over the $u$ variable to reduce the above integral over $K*L$ to an integral over $K\times L$.

We compute that
$$\eqalign{
\int_{u=0}^1 (\pi&-\alpha)A^kB^\ell\,du\cr
&=\int_{u=0}^1 (\pi-\alpha)\left(\vphantom{2^3}\sin\alpha\cos(u(\pi-\alpha))+\cos\alpha
\sin(u(\pi-\alpha))\right)^k\left(\vphantom{2^3}\sin(u(\pi-\alpha))\right)^\ell\,du\cr
&=\int_{\beta=\alpha}^\pi \sin^k(\beta-\alpha)\sin^\ell\beta\ d\beta\cr
&=\ph_{k,\ell}(\alpha),}$$
thanks to the substitution $\beta=\pi-u(\pi-\alpha)$.

This leaves us with 
$$\deg f={-1\over\vol S^n}\int_{K\times L} {\ph_{k,\ell}(\alpha)\over\sin^n\alpha}
\det\left(\x,{\partial\x\over\partial s_1},\ldots,{\partial\x\over\partial s_k},\y,
{\partial\y\over\partial t_1},\ldots,{\partial\y\over\partial t_\ell}\right)
ds_1\cdots ds_k\,dt_1\cdots dt_\ell\,.$$
Recalling from Proposition 3.3 that $\Lk(K^k,L^\ell)=-\deg f$, and introducing the abbreviation 
$$[\x,d\x,\y,d\y]=
\det\left(\x,{\partial\x\over\partial s_1},\ldots,{\partial\x\over\partial s_k},\y,
{\partial\y\over\partial t_1},\ldots,{\partial\y\over\partial t_\ell}\right)
ds_1\cdots ds_k\,dt_1\cdots dt_\ell,$$
we get
$$\Lk(K^k,L^\ell)={1\over \vol S^n}\int_{K\times L} 
{\ph_{k,\ell}(\alpha)\over \sin^n\alpha}\,[\x,d\x,\y,d\y],$$
completing the proof of our theorem.

\bigskip

\noindent\bf \Large{5.\ Proof of the corollary.}\rm
\large
\addtolength{\baselineskip}{2pt}

We begin with the result of our main theorem,
$$\Lk(K^k,L^\ell)={1\over\vol S^n}\int_{K\times L} 
{\ph_{k,\ell}(\alpha)\over \sin^n\alpha}\,[\x,d\x,\y,d\y],\eqno(5.1)$$
where
$$\ph_{k,\ell}(\alpha)=\int_{\beta=\alpha}^\pi \sin^k(\beta-\alpha)\sin^\ell\beta\ d\beta.$$
Putting $-L$ in place of $L$, we get
$$\eqalign{
\Lk(K^k,-L^\ell)&={1\over\vol S^n}\int_{K\times L} 
{\ph_{k,\ell}(\pi-\alpha)\over \sin^n(\pi-\alpha)}\,[\x,d\x,-\y,-d\y]\cr
&={(-1)^{\ell+1}\over\vol S^n}\int_{K\times L} 
{\ph_{k,\ell}(\pi-\alpha)\over \sin^n\alpha}\,[\x,d\x,\y,d\y].}\eqno(5.2)$$
The orientation on $L$ is transferred via the antipodal map of $S^n$ to 
give the orientation on $-L$.

We compute the value of $\ph_{k,\ell}(\pi-\alpha)$ to be
$$\ph_{k,\ell}(\pi-\alpha)=(-1)^k\int_{\beta=0}^\alpha\sin^k(\beta-\alpha)
\sin^\ell\beta\ d\beta.$$
It follows that 
$$\eqalign{
\ph_{k,\ell}(\alpha)+(-1)^k\ph_{k,\ell}(\pi-\alpha)
&=\int_{\beta=0}^\pi\sin^k(\beta-\alpha)\sin^\ell\beta\ d\beta\cr
&=(-1)^k\int_{\beta=0}^\pi \sin^k(\alpha-\beta)\sin^\ell\beta\ d\beta\cr
&=(-1)^k\sin^k*\sin^\ell(\alpha).}\eqno(5.3)$$
Assembling (5.1), (5.2) and (5.3), we get
$$\Lk(K^k,L^\ell)+(-1)^n\Lk(K^k,-L^\ell)={(-1)^k\over\vol S^n}\int_{K\times L}
{\sin^k*\sin^\ell(\alpha)\over\sin^n\alpha}\,[\x,d\x,\y,d\y],$$
completing the proof of the corollary.

\bigskip

\noindent\bf \Large{6.\ Examples}\rm
\large
\addtolength{\baselineskip}{2pt}

\noindent\bf Example 1\it.\ \ Great subspheres $S^k$ and $S^\ell$ in $S^n$\rm.

Our linking integral in this case is
$$\Lk(S^k,S^\ell)={1\over\vol S^n}\int_{S^k\times S^\ell} {\ph_{k,\ell}(\alpha)\over
\sin^n\alpha}\,[\x,d\x,\y,d\y].$$
The geodesic distance $\alpha(\x,\y)$ from each point $\x\in S^k$ to each point $\y\in S^\ell$ is $\pi/2$. Therefore, 
$\sin^n\alpha\equiv 1$, and 
$$\eqalign{
\ph_{k,\ell}(\alpha)=\ph_{k,\ell}(\pi/2)&=
\int_{\beta=\pi/2}^\pi\sin^k(\beta-\pi/2)
\sin^\ell\beta\ d\beta\cr
&=\int_{\theta=0}^{\pi/2}\sin^k\theta\cos^\ell\theta\ d\theta.}$$
Then
$$\eqalign{
\Lk(S^k,S^\ell)&=
{1\over\vol S^n}\left(
\int_{\theta=0}^{\pi/2}\sin^k\theta\cos^\ell
\theta\ d\theta\right)
\int_{S^k\times S^\ell}[\x,d\x,\y,d\y]\cr
&={1\over\vol S^n}\left(
\int_{\theta=0}^{\pi/2}\sin^k\theta\cos^\ell
\theta\ d\theta\right)
\int_{S^k}[\x,d\x]\int_{S^\ell}[\y,d\y]\cr
&={1\over\vol S^n}\left(
\int_{\theta=0}^{\pi/2}\sin^k\theta\cos^\ell
\theta\ d\theta\right)
(\vol S^k)(\vol S^\ell)\cr
&=1.}$$

To go from the first line to the second in this chain of equalities, we used the fact that the $(n+1)\times(n+1)$ determinant $[\x,d\x,\y,d\y]$ consists of a $(k+1)\times(k+1)$ block in the upper left, an $(\ell+1)\times (\ell+1)$ block in the lower right, and zeros elsewhere.

To go from the third line to the fourth, we used the fact that the region between $S^k$ and $S^\ell$ in $S^n$ is filled by hypersurfaces of the form $S^k(\sin\theta)\times S^\ell(\cos\theta)$, whose volumes may be integrated from $\theta=0$ to $\pi/2$ to give the volume of $S^n$.

\medskip

\noindent\bf Example 2\it.\ \ Two curves in $S^3$\rm.

Our linking integral in this case is 
$$\Lk(K^1,L^1)={1\over \vol S^3}\int_{K\times L}{\ph_{1,1}(\alpha)\over\sin^3\alpha}
[\x,d\x,\y,d\y].$$
We compute that 
$$\eqalign{
\ph_{1,1}(\alpha)&=\int_{\beta=\alpha}^\pi\sin(\beta-\alpha)\sin\beta\,d\beta\cr
&={1\over 2}((\pi-\alpha)\cos\alpha+\sin\alpha),}$$
and, using the fact that $\vol S^3=2\pi^2$, get the formula
$$\Lk(K^1,L^1)={1\over 4\pi^2}\int_{K\times L}{(\pi-\alpha)\cos\alpha+\sin\alpha\over
\sin^3\alpha}\,[\x,d\x,\y,d\y],$$
which is equivalent to formula (2) on page 2 of our (2004) paper, and appears as formula 6 on page 7 of Greg Kuperberg's (2006) paper.

\medskip

\noindent\bf Example 3\it.\ \ A curve and a surface in $S^4$\rm.

Our linking integral in this case is
$$\Lk(K^1,L^2)={1\over \vol S^4}\int_{K\times L}
{\ph_{1,2}(\alpha)\over\sin^4\alpha}\,[\x,d\x,\y,d\y].$$
We compute that 
$$\eqalign{
\ph_{1,2}(\alpha)&=\int_{\beta=\alpha}^\pi \sin(\beta-\alpha)\sin^2\beta\ d\beta\cr
&={1\over 3}(1+\cos\alpha)^2,}$$
and, using the fact that $\vol S^4=8\pi^2/3$, get the formula
$$\Lk(K^1,L^2)={1\over 8\pi^2}\int_{K\times L}
{(1+\cos\alpha)^2\over\sin^4\alpha}\,[\x,d\x,\y,d\y].$$

\medskip

\noindent\bf Example 4\it.\ \ Two curves in $S^3$\rm.

We illustrate the corollary in the case that $K^1$ and $L^1$ are closed curves in $S^3$, disjoint from one another, and each from the antipodal image of the other.

The formula from the corollary in this case is
$$\Lk(K^1,L^1)-\Lk(K^1,-L^1) ={-1\over \vol S^3}\int_{K\times L} 
{\sin*\sin(\alpha)\over \sin^3\alpha}\,[\x,d\x,\y,d\y].$$
We compute that 
$$\eqalign{
\sin*\sin(\alpha)&=\int_{\beta=0}^\pi\sin(\alpha-\beta)\sin\beta\ d\beta\cr
&=-{\pi\over 2}\cos\alpha,}$$
and, once again using the fact that $\vol S^3=2\pi^2$, get the formula
$$\Lk(K^1,L^1)-\Lk(K^1,-L^1)={1\over 4\pi}\int_{K\times L}
{\cos\alpha\over\sin^3\alpha}\,[\x,d\x,\y,d\y].$$
When $K^1$ and $L^1$ are orthogonal great circles in $S^3$, then $-L^1=L^1$ with the same orientation, so the left side of the formula directly above is zero. At the same time, the integrand on the right side is identically zero, since $\cos\alpha=\cos \pi/2 = 0$.

\medskip

\noindent\bf Example 5\it.\ \ Two surfaces in $S^5$\rm.

We illustrate the corollary once more, this time in the case that $K^2$ and $L^2$ are closed surfaces in $S^5$, disjoint from one another, and each from the antipodal image of the other.

The formula from the corollary in this case is
$$\Lk(K^2,L^2)-\Lk(K^2,-L^2) ={1\over \vol S^5}\int_{K\times L} 
{\sin^2*\sin^2(\alpha)\over \sin^5\alpha}\,[\x,d\x,\y,d\y].$$
We compute that 
$$\eqalign{
\sin^2*\sin^2(\alpha)&=\int_{\beta=0}^\pi\sin^2(\alpha-\beta)\sin^2\beta\ d\beta\cr
&={\pi\over 8}(1+2\cos^2\alpha),}$$
and, using the fact that $\vol S^5=\pi^3$, get the formula
$$\Lk(K^2,L^2)-\Lk(K^2,-L^2)={1\over 8\pi^2}\int_{K\times L}{1+2\cos^2\alpha\over
\sin^5\alpha}\,[\x,d\x,\y,d\y].$$
When $K^2$ and $L^2$ are orthogonal great 2-spheres in $S^5$ with linking number 1, then $-L^2=L^2$ but has the opposite orientation, so the left side of the formula directly above has the value 2. To check the right side, we have $\alpha=\pi/2$,
so the fraction in the integrand is identically 1. Therefore 
the right side equals
$${1\over 8\pi^2}(\area S^2)(\area S^2)={1\over 8\pi^2}(4\pi)(4\pi) = 2,$$
a reassuring consistency check.

\vfill
\eject

\noindent\bf \Large{References}\rm

\large

\begin{description}
\item[1820] Jean-Baptiste Biot and Felix Savart, {\it Note sur le magnetisme de la pile de Volta},
Annales de chimie et de physique, 2nd ser.,  \bf 15\rm, 222--223.

\item[1833] Carl Friedrich Gauss, {\it Integral formula for linking number}, in
{\it Zur mathematischen theorie der electrodynamische wirkungen},
Collected Works, Vol 5, K\"oniglichen Gesellschaft des Wissenschaften,
Gottingen, 2nd edition, page 605.

\item[1998] Moritz Epple, {\it Orbits of asteroids, a braid, and the first link invariant},
The Mathematical Intelligencer, \bf 20\rm(1)

\item[2004] Dennis DeTurck and Herman Gluck, {\it The Gauss Linking Integral on the 3-sphere
and in hyperbolic 3-space}, announcement, arXiv:math.GT/0406276.

\item[2005] Dennis DeTurck and Herman Gluck, {\it Electrodynamics and the Gauss Linking Integral on the 3-sphere
and in hyperbolic 3-space}, arXiv:math.GT/0510388, with a shortened version to appear in 
J.\ Math.\ Phys.

\item[2006] Greg Kuperberg, {\it From the Mahler conjecture to Gauss linking forms},\\
arXiv:math.MG/0610904v2.

\item[2008] Clayton Shonkwiler and David Shea Vela-Vick, {\it Higher-dimensional linking integrals}, arXiv:0801.4022 [math.GT].

\end{description}

\noindent University of Pennsylvania\\
Philadelphia, PA  19104

\noindent
\it deturck@math.upenn.edu\rm\\
\it gluck@math.upenn.edu\rm\\

\end{document}        

\end{document}